\documentclass[12pt]{article}
\usepackage{amssymb}
\usepackage{latexsym,bm}
\usepackage{graphicx}
\usepackage{amsmath}
\usepackage{mathrsfs}
\usepackage{epstopdf}
\usepackage{tikz}

\date{}
\begin{document}
\title{Every $2$-connected $[4, 2]$-graph of order at least seven contains a pancyclic edge}
\author{\hskip -10mm Chengli Li and  Xingzhi Zhan\thanks{Corresponding author.}}
\maketitle
\footnotetext[1]{Department of Mathematics,  Key Laboratory of MEA (Ministry of Education)
 \& Shanghai Key Laboratory of PMMP, East China Normal University, Shanghai 200241, China}
\footnotetext[2]{E-mail addresses: {lichengli0130@126.com (C. Li),}  {zhan@math.ecnu.edu.cn (X. Zhan).}}

\begin{abstract}
A graph $G$ is called an $[s,t]$-graph if any induced subgraph of $G$ of order $s$ has size at least $t.$ An edge $e$ in a graph $G$ of order $n$ is called
pancyclic if for every integer $k$ with $3\le k\le n,$ $e$ lies in a $k$-cycle. We prove that every $2$-connected $[4, 2]$-graph of order at least seven contains a pancyclic edge.
This strengthens an existing result. We also determine the minimum size of a $[4, 2]$-graph of a given order and show that any $[4, 2]$-graph of order at least eight is not uniquely hamiltonian.
\end{abstract}


\vskip 8mm

\section{Introduction}

We consider finite simple graphs and use standard terminology and notations from [2] and [11]. The {\it order} of a graph is its number of vertices, and the
{\it size} its number of edges.  A $k$-cycle is a cycle of length $k.$ In 1971 Bondy [1] introduced the concept of a pancyclic graph. A graph $G$ of order $n$ is called
{\it pancyclic} if for every integer $k$ with $3\le k\le n,$ $G$ contains a $k$-cycle. See the book [6] for this topic.

{\bf Definition 1.} An edge $e$ of a graph $G$ of order $n$ is said to be {\it pancyclic} if for every integer $k$ with $3\le k\le n,$ $e$ lies in a $k$-cycle of $G.$
A graph $H$ is called {\it edge-pancyclic} if every edge of $H$ is a pancyclic edge.

Clearly, if a graph contains a pancyclic edge, then it is pancyclic. It seems that the concept of a {\it pancyclic edge} is new, although it is very natural.

{\bf Definition 2.} A graph $G$ of order $n$ is said to be {\it vertex-pancyclic} if for every vertex $v$ of $G$ and for every integer $k$ with $3\le k\le n,$ $v$ lies in a $k$-cycle of $G.$

Every edge-pancyclic graph is vertex-pancyclic. The converse is not true. In fact, there exist vertex-pancyclic graphs that contain no pancyclic edge; see an example in Section 3.

We will take the first step in searching for sufficient conditions for a graph to contain a pancyclic edge. It seems that all existing sufficient conditions
for a graph to be pancyclic also ensure that the graph contains a pancyclic edge with possibly very few exceptional graphs.

We consider a class of graphs which have not been fully investigated.

{\bf Definition 3.} Let $s$ be a positive integer and let $t$ be a nonnegative integer. A graph $G$ of order at least $s$ is called an {\it $[s,t]$-graph} if any induced subgraph of $G$ of order $s$ has size at least $t.$

Denote by $\alpha (G)$ the independence number of a graph $G$. We have two facts. (1) Every $[s,t]$-graph is an $[s+1,t+1]$-graph; (2) $\alpha(G)\le k$ if and only if
$G$ is a $[k+1,1]$-graph. Thus the concept of an $[s,t]$-graph is an extension of the independence number.

In 2005 Liu and Wang [8] proved the following result.

{\bf Theorem 1.} {\it Every $2$-connected $[4,2]$-graph of order at least six is hamiltonian.}

There are $2$-connected $[4,2]$-graphs which do not satisfy the Chv\'{a}tal-Erd\H{o}s condition for hamiltonian graphs [12].

In 2008 Liu [9] strengthened Theorem 1 as follows.

{\bf Theorem 2.} {\it  Every $2$-connected $[4,2]$-graph of order at least seven is pancyclic.}

Another proof of Theorem 2 is given in [12].

In this paper we further strengthen Theorem 2 by proving that every $2$-connected $[4, 2]$-graph of order at least seven contains a pancyclic edge.
We also determine the minimum size of a $[4, 2]$-graph of a given order and show that any $[4, 2]$-graph of order at least eight is not uniquely hamiltonian.

In Section 2 we prove the main results, and in Section 3 we mention several related unsolved problems.

For graphs we will use equality up to isomorphism, so $G=H$ means that $G$ and $H$ are isomorphic.

\section{Main results}

We denote by $V(G)$ and $E(G)$ the vertex set and edge set of a graph $G,$ respectively, and denote by $|G|$ and $e(G)$ the order and size of $G,$ respectively.
Thus $|G|=|V(G)|$ and $e(G)=|E(G)|.$ For vertex subsets $S,\, T\subseteq V(G),$ we use $G[S]$ to denote the subgraph of $G$ induced by $S,$
and use $[S, \,T]$ to denote the set of edges having one endpoint in $S$ and the other in $T.$  The neighborhood of a vertex $x$  in $G$ is denoted by $N_G(x).$ The degree of $x$ is denoted by ${\rm deg}_G(x).$ The minimum degree of a graph $G$ is denoted by $\delta(G).$
For $S\subseteq V(G),$ $N_S(x)\triangleq N_G(x)\cap S$ and the degree of $x$ in $S$ is ${\rm deg}_S(x)\triangleq |N_S(x)|.$ We denote by $C_n,$ $P_n$ and $K_n$  the cycle of order $n,$ the path of order $n$ and the complete graph of order $n,$ respectively. Let $K_{s,\, t}$ denote the complete bipartite graph on $s$ and $t$ vertices. $\overline{G}$ denotes the complement of a graph $G.$ The distance between two vertices $u$ and $v$ in a graph $G$ is denoted by $d_G(u, v).$

We denote by $G+H$ the vertex-disjoint union of two graphs $G$ and $H.$ The {\it join} of $G$ and $H,$
denoted by  $G\vee H,$ is the graph obtained from $G+H$ by adding edges joining every vertex of $G$ to every vertex of $H.$

Suppose that $H$ and $R$ are two subgraphs of a graph $G.$ The neighborhood of $H$ in $R,$  denoted by $N_R(H),$ is defined as
$N_R(H)=\{x\in V(R)|\, x\,\, {\rm has}\,\, {\rm a}\,\, {\rm neighbor}\,\, {\rm in}\,\, H\}.$  Given a fixed graph $M,$ we say that a graph $G$ is $M$-free if $G$ does not contain $M$ as a subgraph.

Let $C$ be an oriented cycle, where the orientation is always clockwise. For $u\in V(C),$ denote by $u^{+1}$ the immediate successor of $u$ and $u^{-1}$ the immediate predecessor of $u$ on $C.$ For an integer $\ell\ge 2,$ denote by $u^{+\ell}$ the immediate successor of $u^{+(\ell-1)}$ and $u^{-\ell}$ the immediate predecessor of $u^{-(\ell-1)}$ on $C.$ For convenience, we write $u^+$ for $u^{+1}$ and $u^-$ for $u^{-1}.$ For $x,y\in V(C),$ $x\overrightarrow{C}y$ denotes the path of $C$ from $x$ to $y$ which follows the orientation of $C,$ while $x\overleftarrow{C}y$ denotes the path of $C$ from $x$ to $y$ in the opposite orientation. In particular, if $x=y,$ then $x\overrightarrow{C}y=x\overleftarrow{C}y=x.$

We prepare for the proof of the first main result Theorem 5 below.  Given a graph $H$ and a positive integer $k,$ the {\it $k$-blow-up of $H,$} denoted by $H^{(k)},$ is the graph obtained by replacing every vertex of $H$ with $k$ pairwise nonadjacent vertices where a copy of $u$ is adjacent to a copy of $v$ in the blow-up graph if and only if $u$ is adjacent to $v$ in $H.$

{\bf Lemma 3} [12, Lemma 2.4]. {\it Let $G$ be a $[p+2,p]$-graph of order $n$ with $\delta(G)\ge p\ge 2$ and $n\ge 2p+3.$ Then $G$ is triangle-free if and only if
$p$ is even, $p\ge 6$ and $G=C_5^{(p/2)}.$}

{\bf Definition 4.} Let $G$ be a graph of order $n$ and let $A\subseteq\{3,4,\dots,n\}.$ An edge $e$ of $G$ is called {\it $A$-cyclic} if for every $k\in A,$ $e$ lies
in a $k$-cycle of $G.$

{\bf Notation.} For a graph $G$ and $S\subseteq V(G),$ we denote $e(S)\triangleq e(G[S]).$

{\bf Lemma 4.} {\it Every $2$-connected $[4,2]$-graph of order at least seven contains a $\{3,4,5\}$-cyclic edge.}

{\bf Proof.}
The diamond graph is obtained from the complete graph $K_4$ by removing one edge, and the house graph is obtained from the cycle $C_5$ by adding one chord, which are depicted in Figure 1.
\begin{figure}[h]
\centering
\includegraphics[width=0.55\textwidth]{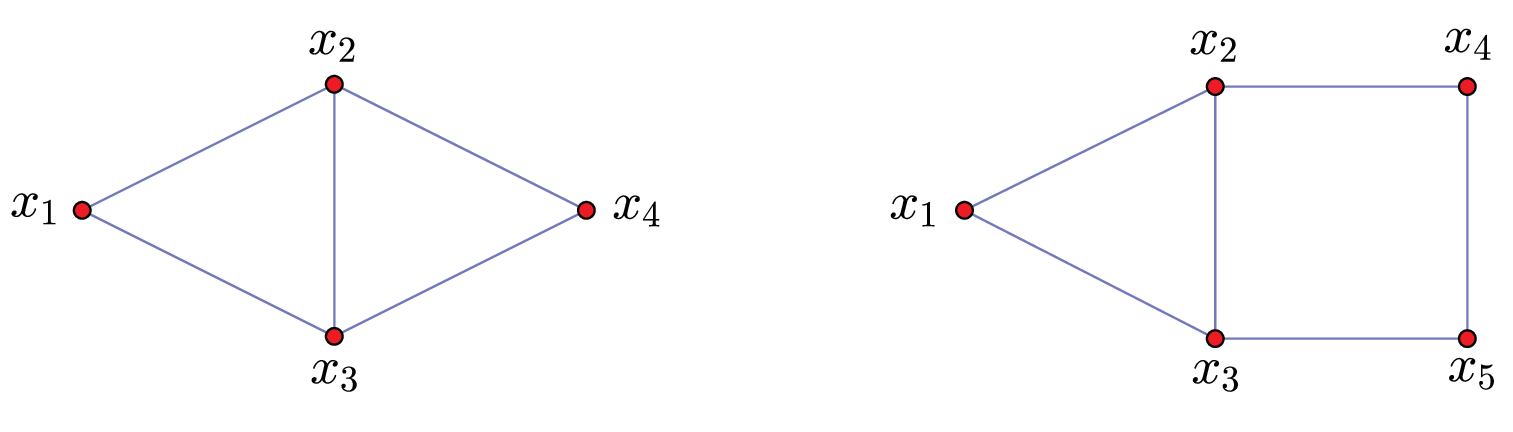}
\caption{The diamond and house}
\end{figure}

Let $G$ be a $2$-connected $[4,2]$-graph of order at least seven. We first prove that $G$ contains either a diamond or a house as a subgraph. By Lemma 3, $G$ contains a triangle, which we denote by $x_1x_2x_3.$ Since $G$ is 2-connected, there exist two nonadjacent edges in $[\{x_1,x_2,x_3\},\,V(G)\setminus \{x_1,x_2,x_3\}].$
Without loss of generality, we assume that these two edges are $x_1y_1$ and $x_2y_2,$ where $y_1\neq y_2$ and $y_1,y_2\in V(G)\setminus \{x_1,x_2,x_3\}.$

Suppose $\deg_G(x_3)\ge 3.$ Let $y_3\in N_G(x_3)\setminus \{x_1,x_2\}.$ If $y_3=y_1$ or $y_3=y_2,$ then $G$ clearly contains a diamond.
Now consider the case when $y_3\neq y_1$ and $y_3\neq y_2.$ By the definition of a $[4,2]$-graph, $e(\{x_1,y_1,y_2,y_3\})\ge 2,$
so at least one of $y_1y_2, y_1y_3,y_2y_3, y_2x_1,y_3x_1$ belongs to $E(G).$ In each of these cases, $G$ contains either a diamond or a house, as desired.

Suppose $\deg_G(x_3)=2.$ If $\{y_1,y_2,x_3\}$ is not an independent set, then $G[\{x_1,x_2,x_3,y_1,y_2\}]$ contains either a house or a diamond.
Next suppose $\{y_1,y_2,x_3\}$ is an independent set. Let $z$ be an arbitrary vertex in $V(G)\setminus \{x_1,x_2,x_3,y_1,y_2\}.$
Since $N_G(x_3)=\{x_1, x_2\}$ and $e(\{x_3,y_1,y_2,z\})\ge 2,$ $z$ must be adjacent to both $y_1$ and $y_2.$

Since $|G|\ge 7,$ there exist two vertices $z_1,z_2$ in $V(G)\setminus \{x_1,x_2,x_3,y_1,y_2\}.$
If $z_1z_2\in E(G),$ then $G[\{z_1,z_2,y_1,y_2\}]$ contains a diamond.
Suppose instead that $z_1z_2\notin E(G).$ Since $e(\{x_1,x_2,z_1,z_2\})\ge 2,$ we deduce that one of
$z_1x_1, z_1x_2,$ $z_2x_1, z_2x_2$ belongs to $E(G).$ In each of these cases, $G$ contains a house, as desired.

Thus, we have proved that $G$ contains either a diamond or a house. Next we will prove Lemma 4 by considering these two cases separately.

{\bf Case 1.} $G$ contains a diamond.

Let $A=\{x_1,x_2,x_3,x_4\}\subseteq V(G)$ such that $G[A]$ contains a diamond, which is depicted in Figure 1.
Since $G$ is 2-connected, there exists a set $\Gamma$ consisting of two nonadjacent edges in $[A,\,V(G)\setminus A].$ By symmetry, we need only consider the following three cases.

{\bf Subcase 1.1.} $\Gamma=\{x_3y_3,x_4y_4\}.$

Since $e(\{x_1,x_2,y_3,y_4\})\ge 2,$  we deduce that one of $y_3y_4, y_3x_1, y_4x_1, y_4x_2, y_3x_2$ belongs to $E(G).$
If one of $y_3y_4,y_3x_1,y_4x_1$ belongs to $E(G),$ then $x_2x_4$ is a $\{3,4,5\}$-cyclic edge. If $y_4x_2\in E(G),$ then $x_1x_2$ is a $\{3,4,5\}$-cyclic edge.
Suppose $y_3x_2\in E(G).$ Since $e(\{x_1,x_4,y_3,y_4\})\ge 2,$ we deduce that one of $y_3x_4,y_3x_1,x_1x_4,x_1y_4,y_3y_4$ belongs to $E(G).$ In each of these cases, $x_2x_4$ is a $\{3,4,5\}$-cyclic edge, as desired.

{\bf Subcase 1.2.} $\Gamma=\{x_2y_2,x_3y_3\}.$

Since $e(\{x_1,x_4,y_2,y_3\})\ge 2,$  we deduce that one of $y_2y_3, x_1y_2,$ $x_4y_2,x_1y_3,$ $x_4y_3$ belongs to $E(G).$
If one of $x_1y_2,x_4y_2,y_2y_3$ belongs to $E(G),$ then $x_1x_3$ is a $\{3,4,5\}$-cyclic edge.
If $x_1y_3\in E(G)$ or $x_4y_3\in E(G),$ then $x_1x_2$ is a $\{3,4,5\}$-cyclic edge.

{\bf Subcase 1.3.} $\Gamma=\{x_1y_1,x_4y_4\}.$

Since $e(\{x_2,x_3,y_1,y_4\})\ge 2,$  we deduce that one of $y_1x_2,y_1x_3,$ $y_4x_2,y_4x_3,$ $y_1y_4$ belongs to $E(G).$
If one of $x_2y_1,x_2y_4,y_1y_4$ belongs to $E(G),$ then $x_1x_3$ is a $\{3,4,5\}$-cyclic edge.
If $x_3y_1\in E(G)$ or $x_3y_4\in E(G),$ then $x_1x_2$ is a $\{3,4,5\}$-cyclic edge.

{\bf Case 2.} $G$ is diamond-free and $G$ contains a house.

Let $B=\{x_1,x_2,x_3,x_4,x_5\}\subseteq V(G)$ such that $G[B]$ contains a house, which is depicted in Figure 1.

{\bf Subcase 2.1.} There exists a vertex $y\in V(G)\setminus B$ such that $\deg_B(y)\ge 2.$

Since $G$ is diamond-free, we need only consider the following cases by symmetry.
If $yx_4,yx_5\in E(G),$ then $x_4x_5$ is a $\{3,4,5\}$-cyclic edge. If $yx_2,yx_4\in E(G),$ then $x_2x_4$ is a $\{3,4,5\}$-cyclic edge. If $yx_1,yx_4\in E(G),$ then $x_1x_2$ is a $\{3,4,5\}$-cyclic edge.
Suppose $yx_3,yx_4\in E(G).$ Since $e(\{x_1,x_2,x_5,y\})\ge 2,$ we deduce that one of  $x_1x_5,x_2x_5,x_1y,x_2y,x_5y$ belongs to $E(G).$ In each of these cases, $G$ contains a diamond, a contradiction.

{\bf Subcase 2.2.} For each vertex $y\in V(G)\setminus B,$ we have $\deg_B(y)\le 1.$

Since $G$ is diamond-free, $x_1x_4,x_1x_5,x_2x_5,x_3x_4 \notin E(G).$ It follows that $G[B]$ is a house.

Suppose that $x_5$ has a neighbor $y_5$ with $y_5\neq x_1,x_2,x_3,x_4.$
Since $\deg_B(y_5)\le 1,$ we have $y_5x_1,y_5x_3,y_5x_4\notin E(G).$ By our assumption, we have $x_1x_4,x_3x_4\notin E(G).$ Thus $G[\{x_1,x_3,x_4,y_5\}]$ contains exactly one edge $x_1x_3,$ a contradiction.
The case that $x_4$ has a neighbor $y_4$ with $y_4\neq x_1,x_2,x_3,x_5$ can be treated by a similar argument.

Suppose that $x_2$ has a neighbor $y_2$ with $y_2\neq x_1,x_3,x_4,x_5.$
Since $\deg_B(y_2)\le 1,$ we have $y_2x_1,y_2x_3,y_2x_4\notin E(G).$ By our assumption, we have $x_1x_4,x_3x_4\notin E(G).$ Thus $G[\{x_1,x_3,x_4,y_2\}]$ contains exactly one edge $x_1x_3,$ a contradiction.
The case that $x_3$ has a neighbor $y_3$ with $y_3\neq x_1,x_2,x_4,x_5$ can be treated by a similar argument.

Among $x_1,x_2,x_3,x_4,x_5,$ only $x_1$ has a neighbor outside $B.$ Hence $x_1$ is a cut vertex of $G,$ which contradicts our assumption that
$G$ is $2$-connected. This completes the proof of Lemma 4.\hfill $\Box$

{\bf Theorem 5.} {\it Every $2$-connected $[4,2]$-graph of order at least seven contains a pancyclic edge.}

{\bf Proof.}
To the contrary, let $G$ be a $2$-connected $[4,2]$-graph of order at least seven that does not contain a pancyclic edge.

Choose an edge $e=uv$ such that $e$ is a $\{3,4,\dots,\ell\}$-cyclic edge and $\ell$ is as large as possible. By Lemma 4 and our assumption, we have $5\le \ell\le n-1.$ By the choice of $e,$ $e$ does not lie in an $(\ell+1)$-cycle. Let $C$ be an $\ell$-cycle containing $e.$ For convenience, we specify an orientation for the cycle $C.$

{\bf Claim 1.} For any vertex $x\in V(G)\setminus V(C),$ $|N_C(x)\setminus \{u,v\}|\le 1.$

To the contrary, we assume that $w$ is a vertex in $V(G)\setminus V(C)$ with two neighbors $v_1,v_2$ such that $\{v_1,v_2\}\cap \{u,v\}=\emptyset.$
Without loss of generality, we may assume that $e$ lies in $v_2^+\overrightarrow{C}v_1^-.$
Clearly, if $v_2=v_1^+,$ then $e$ lies in an $(\ell+1)$-cycle $v_1wv_2\overrightarrow{C}v_1,$ a contradiction.
We will prove Claim 1 by dividing it into the following two cases.

{\bf Case 1.} $|v_1^+\overrightarrow{C}v_2^-|=1.$

In this case, $v_1^+=v_2^-.$ Since $e(\{v_1^-,v_1^+,v_2^+,w\})\ge 2,$ we deduce that one of  $wv_1^+,wv_2^+,wv_1^-,$ $v_2^+v_1^+,v_1^-v_1^+$ belongs to $E(G).$ In each of these cases, $e$ lies in an $(\ell+1)$-cycle, a contradiction.

{\bf Case 2.} $|v_1^+\overrightarrow{C}v_2^-|\ge 2.$

If one of $v_2^-v_1^-,wv_1^-,wv_2^-,wv_2^+$ belongs to $E(G),$ then $e$ lies in an $(\ell+1)$-cycle, a contradiction. Since $e(\{v_1^-,v_2^-,v_2^+,w\})\ge 2,$ we have $v_2^+v_2^-\in E(G).$
If one of $v_2^{-2}v_1^-,wv_1^-,wv_2^{-2},wv_2^+$ belongs to $E(G),$ then $e$ also lies in an $(\ell+1)$-cycle, a contradiction.
Since $e(\{v_1^-,v_2^{-2},v_2^+,w\})\ge 2,$ we have $v_2^+v_2^{-2}\in E(G).$
Repeating the above argument finitely many times, we obtain $v_2^+v_1^+\in E(G),$ and hence $e$ lies in an $(\ell+1)$-cycle $v_1wv_2\overleftarrow{C}v_1^+v_2^+\overrightarrow{C}v_1,$ a contradiction.

This completes the proof of Claim 1.

We say that a component of a graph is nontrivial if its order is at least two. We will use the case $k=2$ of the following fact [5, Lemma 8]:

Let $S$ be a vertex cut of a $k$-connected graph $G$ and let $H$ be a component of $G-S.$ If $|H|\ge k,$ then the set $[S,\, V(H)]$
contains a matching of cardinality $k.$

Let $H$ be a nontrivial component of $G-V(C).$
Since $G$ is $2$-connected, there exist two nonadjacent edges in $[V(C),\, V(H)].$ Moreover, we have the following claim.

{\bf Claim 2.} For any two nonadjacent edges $u_1v_1,u_2v_2$ in $[V(C),\, V(H)]$ where $v_1,v_2\in V(C),$ we have $v_1v_2\in E(C).$

To the contrary, we assume that $v_1v_2\notin E(C).$  Without loss of generality, we assume that $e$ lies in $v_2^+\overrightarrow{C}v_1.$
We distinguish two cases by considering the length of the path $v_1\overrightarrow{C}v_2.$

{\bf Case 1.} $v_2^-=v_1^+.$

In this case, we have $u_1u_2\notin E(G),$ since otherwise $e$ lies in an $(\ell+1)$-cycle $v_1u_1u_2v_2\overrightarrow{C}v_1,$ a contradiction. We will
divide Case 1 into two subcases by considering the position of the edge $e$ on $C.$

{\bf Subcase 1.1.} $e=v_1v_1^-.$

If $u_1v_1^+\in E(G),$ then $e$ lies in an $(\ell+1)$-cycle $v_1u_1v_1^+\overrightarrow{C}v_1,$ a contradiction. Since $\ell\ge 5,$ we have $v_2^+\neq v_1^-.$
By Claim 1, we have $u_2v_1^+,u_2v_2^+\notin E(G).$ Since $e(\{v_2^+,v_1^+,u_1,u_2\})\ge 2,$ we have $v_1^+v_2^+,u_1v_2^+\in E(G).$
One readily observes that $u_1v_2^{+2}\notin E(G).$ If $u_2v_2^{+2}\in E(G),$ then $e$ lies in an $(\ell+1)$-cycle $v_2^+u_1v_1\overleftarrow{C}v_2^{+2}u_2v_2v_2^+,$ a contradiction.
So $G[\{v_2^{+2},v_1^+,u_1,u_2\}]$ contains only one possible edge $v_1^+v_2^{+2},$ a contradiction.

{\bf Subcase 1.2.} $e\neq v_1v_1^-.$

We first consider the case where $\ell=5.$
In this case, $e=v_2^+v_1^-.$
If one of  $u_2v_1^+,u_2v_2^+,u_1v_1^+$ belongs to $E(G),$ then $e$ lies in a $6$-cycle, a contradiction.
Since $e(\{u_1,u_2,v_1^+,v_2^+\})\ge 2,$ we have $u_1v_2^+,v_1^+v_2^+\in E(G).$
Similarly, we also have $u_2v_1^-,v_1^+v_1^-\in E(G).$
Thus $v_2v_2^+$ is a $\{3,4,5,6\}$-cyclic edge, a contradiction.

We next consider the case where $\ell\ge 6.$
In this case, $e$ lies in $v_2^+\overrightarrow{C}v_1^-.$ Without loss of generality, we may assume that $v_2^+\notin \{u,v\}.$
By Claim 1, we have $u_1v_1^+,u_1v_2^+,u_2v_1^+,u_2v_2^{+}\notin E(G).$
So $G[\{v_2^{+},v_1^+,u_1,u_2\}]$ contains only one possible edge $v_1^+v_2^{+},$ a contradiction.

{\bf Case 2.} $v_2^-\neq v_1^+.$

In this case, we have $|v_1^+\overrightarrow{C}v_2^-|\ge 2.$ We also treat this case by considering the position of the edge $e$ on $C.$

{\bf Subcase 2.1.} $e=v_1v_1^-.$

{\bf Subcase 2.1.1.} $v_2^+\neq v_1^-.$

In this case, $|v_2^+\overrightarrow{C}v_1^-|\ge 2.$
We first consider the case where $v_2^{-2}\neq v_1^+.$ This implies that $|v_1^+\overrightarrow{C}v_2^-|\ge 3.$

Suppose $u_1v_1^{+2}\in E(G).$ If $v_1^+v_1^{+3}\in E(G),$ then $e$ lies in an $(\ell+1)$-cycle $v_1u_1v_1^{+2}v_1^+$ $v_1^{+3} \overrightarrow{C}v_1,$ a contradiction.
By Claim 1, we have $u_1v_1^+,u_1v_1^{+3},u_2v_1^+,u_2v_1^{+3}\notin E(G).$ Thus $G[\{v_1^{+3},v_1^+,u_1,u_2\}]$ contains only one possible edge $u_1u_2,$ a contradiction.

Suppose $u_1v_1^{+2}\notin E(G).$
If $v_2^+u_1\notin E(G),$ since $e(\{v_2^+,v_1^{+2},u_1,u_2\})\ge 2,$ we have $u_1u_2,v_1^{+2}v_2^+\in E(G)$
and hence $e$ lies in an $(\ell+1)$-cycle $v_1u_1u_2v_2\overleftarrow{C}v_1^{+2}v_2^+\overrightarrow{C}v_1,$ a contradiction. Thus $v_2^+u_1\in E(G).$
By Claim 1, we have $u_1v_1^{+},u_1v_1^{+3},$ $u_2v_1^+,u_2v_1^{+3}\notin E(G).$
Since $e(\{u_1,u_2,v_1^+,v_1^{+3}\})\ge 2,$ we have $u_1u_2,v_1^+v_1^{+3}\in E(G).$ Thus $e$ lies in an $(\ell+1)$-cycle $v_1^+v_1^{+3}\overrightarrow{C}v_2u_2u_1v_2^+\overrightarrow{C}v_1^+,$ a contradiction.

We next consider the case where $v_2^{-2}= v_1^+.$ This implies that $v_1^+v_2^-\in E(C).$

Suppose $u_1u_2\notin E(G).$
If $u_1v_1^+\in E(G),$ then $e$ lies in an $(\ell+1)$-cycle $v_1u_1v_1^+\overrightarrow{C}v_1,$ a contradiction. Hence $u_1v_1^+\notin E(G).$
By Claim 1, we have $u_2v_2^-,u_2v_1^+\notin E(G).$ Since $e(\{v_1^+,v_2^-,u_1,u_2\})\ge 2,$ we have $u_1v_2^-\in E(G).$ By Claim 1, $u_1v_2\notin E(G).$
If $v_1^+v_2\in E(G),$ then $e$ lies in an $(\ell+1)$-cycle $u_1v_2^-v_1^+v_2\overrightarrow{C}v_1u_1,$ a contradiction. Hence $v_1^+v_2\notin E(G).$
It follows that $G[\{v_2,v_1^+,u_1,u_2\}]$ contains exactly one edge $u_2v_2,$ a contradiction.

Suppose $u_1u_2\in E(G).$  If one of $v_2^-v_2^+,u_2v_2^-,u_2v_2^+$ belongs to $E(G),$ then $e$ lies in an $(\ell+1)$-cycle, a contradiction. Since $e(\{v_2^-,v_2^+,u_1,u_2\})\ge 2,$ we have either $u_1v_2^-\in E(G)$ or $u_1v_2^+\in E(G).$
If $u_1v_2^-\in E(G),$ then $e$ lies in an $\ell$-cycle $C'=v_1u_1v_2^-\overrightarrow{C}v_1.$ Note that $\{v_2,u_1\}\subseteq N_{C'}(u_2)\setminus \{u, v\},$ and
hence $|N_{C'}(u_2)\setminus \{u, v\}|\ge 2.$
On the other hand, in Claim 1 replacing $C$ by $C^{\prime}$ and letting $x=u_2$, we obtain $|N_{C'}(u_2)\setminus \{u, v\}|\le 1,$ a contradiction.
Thus, we have $u_1v_2^+\in E(G).$ If $v_1v_2^-\in E(G),$ then $e$ lies in an $(\ell+1)$-cycle $v_1v_2^-v_2u_2u_1v_2^+\overrightarrow{C}v_1,$ a contradiction.
Recall that we have observed $u_2v_2^+, u_2v_2^-, v_2^-v_2^+\notin E(G).$
Since $e(\{v_2^-,v_2^+,u_2,v_1\})\ge 2,$ we have $u_2v_1,v_1v_2^+\in E(G).$
If $v_2^+v_1^+\in E(G),$ then $e$ lies in an $(\ell+1)$-cycle $v_1u_2v_2v_2^-v_1^+v_2^+\overrightarrow{C}v_1,$ a contradiction.
By Claim 1, we have $u_2v_1^+\notin E(G).$ Thus $G[\{v_2^+,v_2^-,v_1^+,u_2\}]$ contains exactly one edge $v_2^-v_1^+,$ a contradiction.

{\bf Subcase 2.1.2.} $v_2^{+}= v_1^-.$

We first consider the case where $u_1u_2\notin E(G).$

Suppose $\ell\ge 6.$ Then $|v_1^+\overrightarrow{C}v_2^-|\ge 3.$
We assert that $u_1v_1^{+2}\notin E(G).$ To the contrary, assume $u_1v_1^{+2}\in E(G).$
By Claim 1,  $u_1v_1^+, u_1v_1^{+3}, u_2v_1^+, u_2v_1^{+3}\notin E(G).$ We also have $v_1^+v_1^{+3}\notin E(G),$ since otherwise $e$ lies in an $(\ell+1)$-cycle $v_1u_1v_1^{+2}v_1^+v_1^{+3}\overrightarrow{C}v_1,$ a contradiction. Then $G[\{u_1,u_2,v_1^+,v_3^+\}]$ does not contain any edge, a contradiction.
Thus $u_1v_1^{+2}\notin E(G).$
By Claim 1, we have $u_2v_1^+,u_2v_1^{+2}\notin E(G).$
We also have $u_1v_1^+\notin E(G),$ since otherwise $e$ lies in an $(\ell+1)$-cycle $v_1u_1v_1^+\overrightarrow{C}v_1,$ a contradiction.
It follows that $G[\{v_1^+,v_1^{+2},u_1,u_2\}]$ contains only one edge $v_1^+v_1^{+2},$ a contradiction.

Suppose $\ell=5.$ Similarly as above, we have $u_1v_1^+\notin E(G).$ By Claim 1, we have $u_2v_1^{+2},u_2v_1^+\notin E(G).$ Since $e(\{v_1^+,v_1^{+2},u_1,u_2\})\ge 2,$ we have $u_1v_1^{+2}\in E(G).$
If $u_1v_2\in E(G),$ then $e$ lies in a $6$-cycle $v_1^{+2}u_1v_2v_2^+v_1v_1^+v_1^{+2},$ a contradiction.
If $v_1^+v_2\in E(G),$ then $e$ lies in a $6$-cycle $v_2^+v_1u_1v_1^{+2}v_1^+v_2v_2^+,$ a contradiction.
Thus $G[\{u_2,v_2,v_1^+,u_1\}]$ contains exactly one edge $u_2v_2,$ a contradiction.

We next consider the case where $u_1u_2\in E(G).$

Suppose $\ell\ge 7.$ Then $|v_1^+\overrightarrow{C}v_2^-|\ge 4.$
We have $v_2^+v_1^{+2}\notin E(G),$ since otherwise $e$ lies in an $(\ell+1)$-cycle $v_2^+v_1^{+2}\overrightarrow{C}v_2u_2u_1v_1v_2^+,$ a contradiction. We also have $u_2v_2^+\notin E(G),$ since otherwise $e$ lies in an $(\ell+1)$-cycle $u_2v_2^+\overrightarrow{C}v_2u_2,$ a contradiction. By Claim 1, we have $u_2v_1^+,u_2v_1^{+2}\notin E(G).$
Since $e(\{u_2,v_2^+,v_1^+,v_1^{+2}\})\ge 2,$ we have $v_2^+v_1^+\in E(G).$
Similarly, since $e(\{u_2,v_2^+,v_1^{+3},v_1^{+2}\})\ge 2,$ we also have $v_2^+v_1^{+3}\in E(G).$

If $u_1v_1^{+2}\in E(G),$ then there exists an $\ell$-cycle $C''=v_1u_1v_1^{+2}\overrightarrow{C}v_1$ containing $e.$ Note that $\{v_2,u_1\}\subseteq N_{C''}(u_2)\setminus \{u,v\},$ and hence $|N_{C''}(u_2)\setminus \{u,v\}|\ge 2.$ On the other hand, in Claim 1 replacing $C$ by $C^{''}$ and letting $x=u_2$, we obtain $|N_{C''}(u_2)\setminus \{u, v\}|\le 1,$ a contradiction. Thus $u_1v_1^{+2}\notin E(G).$ One readily observes that $u_2v_2^+,u_2v_1^{+2}\notin E(G).$
Since $e(\{u_1,u_2,v_1^{+2},v_2^+\})\ge 2,$ we have $u_1v_2^+\in E(G).$
This implies that $v_1^{+2}v_1^{+4}\notin E(G),$
since otherwise $e$ lies in an $(\ell+1)$-cycle $v_2^+\overrightarrow{C}v_1^{+2}v_1^{+4}\overrightarrow{C}v_2u_2u_1v_2^+,$ a contradiction.
By a similar argument, we have $u_2v_1^{+2},u_2v_1^{+4},u_2v_2^+\notin E(G).$ Also we have already shown that $v_2^+v_1^{+2}\notin E(G).$
It follows that $G[\{v_1^{+2},v_1^{+4},u_2,v_2^+\}]$ contains only one possible edge $v_2^+v_1^{+4},$ a contradiction.

Suppose $\ell=6.$
By Claim 1, we have $u_2v_1^+,u_2v_1^{+2}\notin E(G).$
If $u_2v_2^+\in E(G)$ or $v_2^+v_1^{+2}\in E(G),$ then we obtain a $7$-cycle containing $e,$ a contradiction.
Since $e(\{v_2^+,v_1^+,v_1^{+2},$ $u_2\})\ge 2,$ we have $v_2^+v_1^+\in E(G).$
Similarly, since $e(\{v_2^+,v_1^{+3},v_1^{+2},u_2\})\ge 2,$ we have $v_2^+v_1^{+3}\in E(G).$
Then it is easy to see that $v_2^+v_1^{+3}$ is a $\{3,4,5,6,7\}$-cyclic edge, a contradiction.

Suppose $\ell=5.$
Using an argument similar to the one above, we have
$$
u_2v_1^+,u_2v_1^{+2},u_2v_2^+,v_2^+v_1^{+2}\notin E(G).
$$
Since $e(\{u_2,v_1^+,v_1^{+2},v_2^+\})\ge 2,$ we have $v_2^+v_1^+\in E(G).$
Since $e(\{u_2,u_1,v_1^{+2},v_2^+\})\ge 2,$ we have either $u_1v_1^{+2}\in E(G)$ or $u_1v_2^+\in E(G).$ In either case, $v_1v_1^+$ is a $\{3,4,5,6\}$-cyclic edge, a contradiction.

{\bf Subcase 2.2.} $e\neq v_1v_1^-.$

In this case, $e$ lies in $v_2^+\overrightarrow{C}v_1^-.$
By Claim 1, we have $u_1v_1^+,u_1v_1^{+2},u_2v_1^+,u_2v_1^{+2}\notin E(G).$ Since $e(\{u_1,u_2,v_1^+,v_1^{+2}\})\ge 2,$ we deduce that $u_1u_2 \in E(G).$
Next we assert that $V(v_1^+\overrightarrow{C}v_2^-)$ is a clique in $G.$ Suppose not. Let $z_1,z_2$ be two nonadjacent vertices in $v_1^+\overrightarrow{C}v_2^-.$ Then $G[\{u_1,u_2,z_1,z_2\}]$ contains exactly one edge $u_1u_2,$ a contradiction.

By symmetry, we need only consider the following two cases.

{\bf Subcase 2.2.1.} $e=v_2^+v_1^-.$

Since $v_2^-\neq v_1^+,$ we have $\ell\ge 6.$
Suppose $\ell\ge 7.$ By Claim 1, we have $u_1v_1^{+2},u_2v_1^{+2}\notin E(G).$ If $u_2v_2^+\in E(G),$ then $e$ lies in an $(\ell+1)$-cycle $v_2u_2v_2^+\overrightarrow{C}v_2,$ a contradiction. If $u_1v_2^+\in E(G),$ then $e$ lies in an $(\ell+1)$-cycle $ v_2^+\overrightarrow{C}v_1^+v_1^{+3}\overrightarrow{C}v_2u_2u_1v_2^+$
where $v_1^+v_1^{+3}\in E(G)$ since $V(v_1^+\overrightarrow{C}v_2^-)$ is a clique in $G,$ a contradiction.
If $v_1^{+2}v_2^+\in E(G),$ then $e$ lies in an $(\ell+1)$-cycle $
v_2^+v_1^-v_1u_1u_2v_2\overleftarrow{C}v_1^{+2}v_2^+,$ a contradiction.
Thus $G[\{v_2^+,v_1^{+2},u_1,u_2\}]$ contains exactly one edge $u_1u_2,$ a contradiction.

Suppose $\ell=6.$
If $u_2v_2^+\in E(G),$ then $e$ lies in a $7$-cycle $v_2^+\overrightarrow{C}v_2u_2v_2^+,$ a contradiction.
If $v_2^+v_2^-\in E(G),$ then $e$ lies in a $7$-cycle $v_2^+v_1^-v_1u_1u_2v_2v_2^-v_2^+,$ a contradiction.
By Claim 1, we have $u_2v_2^-,u_2v_1^+\notin E(G).$
Since $e(\{v_2^-,v_1^{+},v_2^{+},u_2\})\ge 2,$ we have $v_2^+v_1^+\in E(G).$ Similarly, we obtain $v_2^-v_1^-\in E(G).$
By Claim 1, we have $u_2v_2^-,u_2v_1\notin E(G).$ Also, we have just proved that $u_2v_2^+, v_2^+v_2^-\notin E(G).$
Since $e(\{v_2^-,v_2^{+},v_1,u_2\})\ge 2,$ we have $v_1v_2^-\in E(G).$
Thus $e$ lies in a $7$-cycle $v_2^+v_1^-v_2^-v_1u_1u_2v_2v_2^+,$ a contradiction.

{\bf Subcase 2.2.2.} $v_2^+\notin \{u,v\}.$

By Claim 1, we have $v_2^+u_2,v_2^+u_1,v_1^{+2}u_2,v_1^{+2}u_1\notin E(G).$
If $v_2^+v_1^{+2}\in E(G),$ then $e$ lies in an $(\ell+1)$-cycle $v_2^+\overrightarrow{C}v_1u_1u_2v_2\overleftarrow{C}v_1^{+2}v_2^+,$ a contradiction.
Hence we have $v_2^+v_1^{+2}\notin E(G).$ Now $G[\{v_2^+,v_1^{+2},u_1,u_2\}]$ contains exactly one edge $u_1u_2,$ a contradiction.

This completes the proof of Claim 2.

{\bf Claim 3.} Let $H$ be a nontrivial component of $G-V(C).$ Then $G[N_C(H)]$ is an edge of $C.$

Since $G$ is $2$-connected, we have $|N_C(H)|\ge 2.$
Let $u_1v_1,u_2v_2$ be two nonadjacent edges in $[V(C),V(H)],$ where $v_1,v_2\in V(C).$ By Claim 2, we have $v_1v_2\in E(C).$
Renaming the vertices $u_1,v_1,u_2,v_2$ if necessary, we may assume that $v_2=v_1^+.$
Thus, to prove Claim 3 it suffices to show that $N_C(H)=\{v_1,v_2\}.$ To the contrary, we assume that $v_3$ is a vertex in $N_C(H)$ distinct from $v_1$ and $v_2.$

We first assert that $v_3=v_1^-$ or $v_3=v_2^+.$ Otherwise $v_1u_1, v_3u_3$ or $v_2u_2, v_3u_3$ would be two nonadjacent edges where $u_3\in N_H(v_3),$ but $v_1v_3\notin E(C)$ and $v_2v_3\notin E(C),$ which contradicts the conclusion of Claim 2.
Without loss of generality, we assume that $v_3=v_2^+.$
Applying Claim 2 repeatedly, we have $v_3u_1\in E(G), v_1u_2,v_3u_2\notin E(G)$ and $v_1^-\notin N_C(H).$ Thus $N_C(H)=\{v_1,v_2,v_3\}.$

By Claim 1, we have $v_1\in \{u,v\}$ or $v_3\in \{u,v\}.$ Without loss of generality, we assume that $v_1=u.$
Clearly, we have $u_1v_2\notin E(G),$ since otherwise $e$ lies in an $(\ell+1)$-cycle $v_2u_1v_3\overrightarrow{C}v_2,$ a contradiction.
There are two possible positions of the vertex $v$ on the cycle $C.$

{\bf Case 1.} $v=v_1^-.$

We first consider the case where $\ell\ge 6.$ Note that $e$ lies in a new $\ell$-cycle $v_1u_1v_3\overrightarrow{C}v_1.$
Replacing $C$ in Claim 1 by this new cycle, we deduce that  $N_C(v_2)\setminus \{v,u,v_3\}=\emptyset.$  Since $u_1v_3\in E(G),$ Claim 1 also implies that
$u_1v^{-}, u_1v^{-2}\notin E(G).$ Then $G[\{u_1,v_2,v^-,v^{-2}\}]$ contains exactly one edge $v^-v^{-2},$ a contradiction.

Suppose $\ell=5.$ We have $u_1u_2\in E(G),$ since otherwise $G[\{u_1,u_2,v,v^-\}]$ contains exactly one edge $vv^-,$ a contradiction.
Since $e(\{v,v^-,v_2,u_1\})\ge 2,$ we have either $v_2v^-\in E(G)$ or $v_2v\in E(G).$ In the former case, $e$ lies in a $6$-cycle $vuu_1u_2v_2v^-v,$ a contradiction.
In the later case, it is easy to verify that $v_2v$ is a $\{3,4,5,6\}$-cyclic edge, a contradiction.

{\bf Case 2.} $v=v_2.$

We first consider the case where $\ell\ge 6.$
Since $e(\{u^-,u^{-3},u_1,u_2\})\ge 2,$ we have $u_1u_2,u^-u^{-3}\in E(G).$ Thus $e$ lies in an $(\ell+1)$-cycle $vu_2u_1v_3\overrightarrow{C}u^{-3}u^-uv,$ a contradiction.

Suppose $\ell=5.$ We have $u_1u_2\in E(G),$ since otherwise $G[\{u^-,u^{-2},u_1,u_2\}]$ contains exactly one edge $u^-u^{-2},$ a contradiction.
If $v_3u^-\in E(G),$ then $e$ lies in an $6$-cycle $v_3u^-uvu_2u_1v_3,$ a contradiction. Hence $v_3u^-\notin E(G).$
It follows that $v_3u\in E(G),$ since otherwise $G[\{u,u^-,v_3,u_2\}]$ contains exactly one edge $uu^-,$ a contradiction.
In the paragraph preceding Case 1, we have proved that $u_1v_2\notin E(G).$
Since $e(\{u^-,u^{-2},v,u_1\})\ge 2,$ we have either $vu^-\in E(G)$ or $vu^{-2}\in E(G).$ In either case, we can verify that $u_1u$ is a $\{3,4,5,6\}$-cyclic edge, a contradiction.

This completes the proof of Claim 3.

Let $H_1,H_2,\dots,H_t$ be the nontrivial components of $G-V(C)$ and let $T_1,T_2,\dots,T_k$ be the trivial components of $G-V(C).$
Since $G$ is a $[4,2]$-graph, we have $\alpha(G)\le 3$ and hence $t+k\le 3.$

If $t=3,$ then there exists an induced subgraph isomorphic to $P_2+ 2P_1$ in $H_1+ H_2+ H_3,$ a contradiction.
If $t=2,$ then by Claim 3 and the fact that $\ell\ge 5,$ there exists a vertex in $V(G)\setminus (N_C(H_1)\cup N_C(H_2)),$ and hence $G$ contains an induced subgraph isomorphic to $P_2+ 2P_1,$ a contradiction.

If $t=1,$ then $k\le 1,$ since otherwise there exists an induced subgraph isomorphic to $P_2+ 2P_1$ in $H_1+ T_1+ T_2,$ a contradiction. Let $zz'$ and $ww'$ be the two nonadjacent edges in $[V(C),\, V(H_1)],$ where $z,w\in V(C)$ and $z',w'\in V(H_1).$
Assume that $w=z^+.$

Suppose $k=1.$ If $\ell\ge 6,$ then by Claims 1 and 3, there exists a vertex in $V(C)\setminus (N_C(H_1)\cup N_C(T_1))$ and hence $G$ contains an induced subgraph isomorphic to $P_2+ 2P_1,$ a contradiction.
If $\ell=5,$ then $N_C(H_1)\cup N_C(T_1)=V(C),$ since otherwise we also obtain a vertex in $V(C)\setminus (N_C(H_1)\cup N_C(T_1))$ and hence $G$ contains an induced subgraph isomorphic to $P_2+ 2P_1,$ a contradiction.
Moreover, by Claims 1 and 3, $N_C(H_1)$ consists of two consecutive vertices on $C,$ and $N_C(T_1)$ consists of three consecutive vertices on $C.$
Consequently, $e$ lies in a $6$-cycle in $G,$ a contradiction.

Suppose $k=0.$ Then $V(H_1)$ and $V(C)\setminus \{z,w\}$ are cliques, since otherwise, we can obtain an induced subgraph isomorphic to $P_2+ 2P_1$ in $G,$ a contradiction.
If $\ell\ge 6,$ then $w^+w^{+2}$ is a pancyclic edge, a contradiction.
Suppose $\ell=5.$ If $w^+z\in E(G),$ then $ww^+$ is a $\{3,4,5,6\}$-cyclic edge, a contradiction. If $w^{+2}z\in E(G),$ then $w^{+2}w^+$ is a $\{3,4,5,6\}$-cyclic edge, a contradiction.  If $w'z\in E(G),$ then $ww'$ is a $\{3,4,5,6\}$-cyclic edge, a contradiction. If $w'w^+\in E(G),$ then $ww^+$ is a $\{3,4,5,6\}$-cyclic edge, a contradiction.
If $w'w^{+2}\in E(G),$ then $w^{+2}w^+$ is a $\{3,4,5,6\}$-cyclic edge, a contradiction.
Thus $G[\{w^+,w^{+2},z,w'\}]$ contains only one  edge $w^+w^{+2},$ a contradiction.

Next we consider the case where $t=0.$ Then $1\le k\le 3.$
Recall that $e=uv$ is a $\{3,4,\dots,\ell\}$-cyclic edge. We assume that $v=u^+.$

Suppose $k=3.$ Let $z_i$ be the vertex in $T_i$ for $i=1,2,3.$
If $\ell\ge 6,$ then by Claim 1, there exists a vertex in $V(C)\setminus (\bigcup_{i=1}^3N_C(z_i)).$ Hence $G$ contains an independent set with size four, a contradiction. If $\ell=5,$ then $V(C)=\bigcup_{i=1}^3N_C(z_i),$ since otherwise, $G$ also contains an independent set with size four, a contradiction.
Let $z_i'$ be the unique vertex in $N_C(z_i)\setminus \{u,v\}$ for $i=1,2,3.$ Then $z_1',z_2',z_3'$ are three distinct vertices.
Without loss of generality, we assume that $C=uvz_1'z_2'z_3'u.$ If $z_1v\in E(G)$ or $z_3u\in E(G),$ then $e$ lies in a $6$-cycle, a contradiction.
Since $G$ is $2$-connected, by Claim 1, $z_1u,z_3v\in E(G)$ and $z_2$ is adjacent to at least one of $u$ and $v.$ Without loss of generality, we assume that $z_2u\in E(G).$
Since $e(\{z_1,z_2,z_3',v\})\ge 2,$ we have $vz_3',vz_2\in E(G).$ Thus $vz_2$ is a $\{3,4,5,6\}$-cyclic edge, a contradiction.

Suppose $k=2.$ Let $z_i$ be the vertex in $T_i$ for $i=1,2.$
If $\ell\ge 6,$ then there exist two vertices $x, y$ on $C$ such that $\{x,y\}\cap (N_C(z_1)\cup N_C(z_2))=\emptyset.$
By Claim 1, $G[\{x,y,z_1,z_2\}]$ contains at most one possible edge $xy,$ a contradiction.

Suppose $\ell=5.$ Now $|G|=7.$ We assert that $|V(C)\setminus (N_C(z_1)\cup N_C(z_2))|=1.$ Indeed, by Claim 1, we have $|N_C(z_1)\cup N_C(z_2)|\le 4,$ and hence $|V(C)\setminus (N_C(z_1)\cup N_C(z_2))|\ge 1.$ If $|V(C)\setminus (N_C(z_1)\cup N_C(z_2))|\ge 2,$ let $x,y\in V(C)\setminus (N_C(z_1)\cup N_C(z_2)).$  Then $G[\{x,y,z_1,z_2\}]$ contains only one possible edge $xy,$ a contradiction. This shows that $|V(C)\setminus (N_C(z_1)\cup N_C(z_2))|=1.$
Thus we may assume that either $z_1v^+,z_2v^{+2}\in E(G)$ or $z_1v^+,z_2v^{+3}\in E(G).$
Clearly, $z_1v\notin E(G),$ since otherwise $e$ lies in a $6$-cycle $uvz_1v^+v^{+2}v^{+3}u,$ a contradiction.
In the former case, since $e(\{v^{+3},v,z_1,z_2\})\ge 2,$ we have $vv^{+3},vz_2\in E(G).$
Since $G$ is $2$-connected, by Claim 1 and the fact that $z_1v\notin E(G),$ we have $z_1u\in E(G).$ Thus $e$ lies in a $6$-cycle $uvz_2v^{+2}v^+z_1u,$ a contradiction.
In the latter case, we have $z_2u\notin E(G),$ since otherwise $e$ lies in a $6$-cycle $v^{+3}z_2uvv^+v^{+2}v^{+3},$ a contradiction. Since $e(\{v^{+2},u,z_1,z_2\})\ge 2,$ we have $uz_1,uv^{+2}\in E(G).$ By symmetry, we also have $vz_2,vv^{+2}\in E(G).$ Thus $uv^{+3}$ is a pancyclic edge, a contradiction.

Suppose $k=1.$ Let $z_1$ be the vertex in $T_1.$
Since $G$ is $2$-connected, we have $|N_C(z_1)|\ge 2.$
By Claim 1, there are the following two cases.

{\bf Case 1.} $|N_C(z_1)\setminus \{u,v\}|=1.$

Let $w$ be the neighbor of $z_1$ in $V(C)\setminus \{u,v\}.$
Without loss of generality, we assume that $z_1v\in E(G).$
One readily observes that $w\neq v^+,$ since otherwise $e$ lies in an $(\ell+1)$-cycle $vz_1w\overrightarrow{C}v,$ a contradiction.
We assert that $w\neq u^-.$ Suppose not. This implies that $z_1u\notin E(G),$ since otherwise $e$ lies in an $(\ell+1)$-cycle $wz_1u\overrightarrow{C}w,$ a contradiction.
Let $x$ be an arbitrary vertex in $V(G)\setminus \{w,z_1,u,v,v^+\}.$ By Claim 1, we have $z_1x,z_1v^+\notin E(G).$ We also have $uv^+\notin E(G),$ since otherwise $e$ lies in an $(\ell+1)$-cycle $z_1vuv^+\overrightarrow{C}wz_1,$ a contradiction. Since $e(\{z_1,u,v^+,x\})\ge 2,$ we have $xv^+,xu\in E(G).$ It follows that $ww^-$ is a pancyclic edge, a contradiction. Hence $w\neq u^-.$

Let $y$ be an arbitrary vertex in $V(G)\setminus \{z_1,u,v,w,v^+,w^+\}.$ By Claim 1, we have $z_1v^+,z_1w^+,z_1y\notin E(G).$
We also have $v^+w^+\notin E(G)$ since otherwise $e$ lies in an $(\ell+1)$-cycle $w\overleftarrow{C}v^+w^+\overrightarrow{C}vz_1w,$ a contradiction.
Since $e(\{z_1,w^+,v^+,y\})\ge 2,$ we have $yv^+,yw^+\in E(G).$
If $|v^+\overrightarrow{C}w^-|\ge 3,$ then $e$ lies in an $(\ell+1)$-cycle $z_1ww^-v^+\overrightarrow{C}w^{-2}w^+\overrightarrow{C}vz_1,$ a contradiction.
If $|w^+\overrightarrow{C}u^-|\ge 2,$ then $e$ lies in an $(\ell+1)$-cycle $z_1ww^+w^-\overleftarrow{C}v^+w^{+2}\overrightarrow{C}vz_1,$ a contradiction.
Thus $|v^+\overrightarrow{C}w^-|= 2, |w^+\overrightarrow{C}u^-|=1,$ and hence $|G|=7.$
By what we have proved above, $w^-w^+\in E(G).$
If $uz_1\in E(G),$ then we can verify that $uz_1$ is a pancyclic edge, a contradiction.
So  $uz_1\notin E(G).$
Since $e(\{u,v^+,w^+,z_1\})\ge 2,$ we have $uv^+\in E(G).$
Thus $e$ lies in a $7$-cycle $vuv^+w^-w^+wz_1v,$ a contradiction.

{\bf Case 2.} $N_C(z_1)=\{u,v\}.$

Since $uv$ is a $\{3,4,\dots,n-1\}$-cyclic edge, we deduce that $z_1v$ is a $\{3,5,6,\dots,n\}$-cyclic edge.
If $z_1v$ lies in a $4$-cycle, then $z_1v$ is a pancyclic edge, a contradiction. Let $V'=V(G)\setminus \{z_1,v,u\}.$
Then for any vertex $x\in V',$ either $xv\notin E(G)$ or $xu\notin E(G).$
Let $V_1=\{x\in V':xv\notin E(G)\}$ and let $V_2=\{x\in V':xu\notin E(G)\}.$ Clearly, $V(G)=V_1\cup V_2\cup \{u,v,z_1\}.$
Without loss of generality, we assume that $|V_1|\ge |V_2|.$
We assert that $V_1$ and $V_2$ are cliques, since otherwise $G[V_1\cup \{z_1,v\}]$ or $G[V_2\cup \{z_1,u\}]$ contains an induced subgraph isomorphic to $P_2+ 2P_1,$ a contradiction.

{\bf Subcase 2.1.} $V_1\cap V_2=\emptyset.$

In this case, $v$ is adjacent to each vertex of $V_2$ and $u$ is adjacent to each vertex of $V_1.$
Since $G$ is $2$-connected, there exists an edge in $[V_1,\,V_2],$ say $v'u'$ where $v'\in V_2,u'\in V_1.$
If $|V_1|\ge 3,$ let $u''$ be a vertex in $V_1\setminus \{u'\}.$ Then $u'u''$ is a pancyclic edge, a contradiction.
Since $|G|\ge 7,$ we have $|V_1|=|V_2|=2$ and $|G|=7.$ Let $v''$ be the vertex in $V_2\setminus \{v'\}.$
If one of $v''u',v''u'',v''z_1$ belongs to $E(G),$  then $vv''$ is a pancyclic edge, a contradiction.
If one of $u'z_1,u''z_1$ belongs to $E(G),$  then $uu''$ is a pancyclic edge, a contradiction.
Thus $G[\{z_1,v'',u',u''\}]$ contains exactly one edge $u'u'',$ a contradiction.

{\bf Subcase 2.2.} $V_1\cap V_2\neq \emptyset.$

In this case, $v$ is adjacent to each vertex of $V_2\setminus V_1$ and $u$ is adjacent to each vertex of $V_1\setminus V_2.$
We observe that $V_1\setminus V_2\neq \emptyset$ and $V_2\setminus V_1\neq \emptyset,$ since otherwise $v$ or $u$ would be a cut vertex, a contradiction.
Let $w$ be a vertex in $V_1\cap V_2$ and let $w'$ be a vertex in $V_1\setminus V_2.$
Since $|G|\ge 7,$ we have $|V_1\cup V_2|\ge 4$ and $|V_1|\ge 3,$ which can be seen by considering the two cases $|V_1\cap V_2|\ge 2$ and $|V_1\cap V_2|=1$ separately.
Now it is easy to verify that $ww'$ is a pancyclic edge, a contradiction.
This completes the proof of Theorem 5. \hfill $\Box$

{\bf Remark 1}. In view of Theorem 5, one may wonder whether every $2$-connected $[4, 2]$-graph is vertex-pancyclic.
The answer is no. For example, the $2$-connected $[4, 2]$-graph $(2K_1)\vee (K_1+K_{n-3})$ is not vertex-pancyclic.

Now we consider the following problem.

{\bf Problem 1} (X. Zhan, April 2025, see [4, Problem 1]). Determine the minimum size of a connected $[s, t]$-graph of order $n,$ and determine the extremal graphs.

In the next Theorems 9, 10, 11 we solve the case $s=4$ and $t=2$ of Problem 1 with further connectivity conditions.

The following lemma is a well-known theorem of Mantel.

{\bf Lemma 6} [11, p.41]. {\it If $G$ is a triangle-free graph of order $n,$ then $e(G)\le \lfloor n^2/4\rfloor.$ Furthermore, equality holds if and only if $G=K_{\lfloor n/2\rfloor,\lceil n/2\rceil}.$
}

Brandt [3] generalized Mantel's theorem to the case of nonbipartite graphs as follows.

{\bf Lemma 7} [3, Lemma 2]. {\it Every nonbipartite graph of order $n$ with more than $(n-1)^2/4+1$ edges contains a triangle.
}

Let $K_{\lfloor n/2 \rfloor,\lceil n/2 \rceil}$ be a balanced complete bipartite graph with partite sets $V_1=\{u,y,\dots\}$ and  $V_2=\{v,x,\dots\}.$

{\bf Notation.} $$
G_0=K_{\lfloor n/2 \rfloor,\lceil n/2 \rceil}-uv,\quad G_1 = G_0 - ux, \quad G_2 =G_0- vy,\quad G_3 = G_0 - xy.
$$

Denote by $K_{s, t}^-$ the graph obtained from $K_{s, t}$ by deleting one edge.

{\bf Lemma 8.} {\it Let $G$ be a triangle-free graph of order $n.$ If $n\ge 6$ and $e(G)=\lfloor n^2/4\rfloor-1,$ then $G=G_0,$ or $G=K_{ n/2 -1, n/2 +1}$ ($n$ is even).
If $n\ge 8$ and $e(G)= \lfloor n^2/4 \rfloor-2,$ then $G$ is isomorphic to one of $G_1,G_2,G_3,K_{ n/2 -1, n/2 +1}^-$ ($n$ is even), $K_{ (n-3)/2 ,(n+3)/2}$ ($n$ is odd).
}

{\bf Proof.} One readily observes that $\lfloor n^2/4 \rfloor-1>(n-1)^2/4+1$ when $n\ge 6$ and $\lfloor n^2/4\rfloor-2>(n-1)^2/4+1$ when $n\ge 8.$
So by Lemma 7, $G$ is a bipartite graph. Thus, the conclusion in Lemma 8 can be deduced by simple calculations.
\hfill $\Box$

{\bf Theorem 9.} {\it Let $G$ be a $[4,2]$-graph of order $n.$ Then
$$
e(G)\ge \lfloor (n-1)^2/4\rfloor.\eqno (1)
$$
If $n\ge 7,$ equality holds in (1) if and only if $G=K_{\lfloor n/2 \rfloor}+
K_{\lceil n/2 \rceil}.$
}

{\bf Proof.} We use induction on the order of $G.$
It is easy to verify that the conclusion of Theorem 9 holds for $4\le n\le 6.$

Next let $G$ have order $n\ge 7$ and assume that Theorem 9 holds for all graphs of order less than $n.$
Since $G$ is a $[4,2]$-graph, we have $\alpha(G)\le 3.$
If $\alpha(G)=1,$ then $G$ is a complete graph and hence $e(G)=n(n-1)/2>\lfloor (n-1)^2/4 \rfloor.$

If $\alpha(G)=2,$ then $\overline{G}$ is triangle-free.
By Lemma 6, we have
\[
e(G)=\frac{n(n-1)}{2}-e(\overline{G})\ge \frac{n(n-1)}{2}-\left\lfloor\frac{n^2}{4}\right\rfloor=\left\lfloor\frac{(n-1)^2}{4}\right\rfloor,
\]
with equality if and only if $\overline{G}=K_{\lfloor n/2\rfloor,\lceil n/2\rceil},$ i.e., $G=K_{\lfloor n/2\rfloor}+
K_{\lceil n/2\rceil}.$

Suppose $\alpha(G)=3.$
Let $A$ be an independent set of cardinality three. Since $G$ is a $[4,2]$-graph,
for every $v\in V(G)\setminus A,$ we have $|N_A(v)|\ge 2.$ By the induction hypothesis and the assumption that $n\ge 7,$ we have
\begin{align*}
e(G)\ge e(G-A)+2(n-3)
&\ge \left\lfloor \frac{(n-4)^2}{4}\right\rfloor+2(n-3)\\
&>\left\lfloor\frac{(n-1)^2}{4}\right\rfloor.
\end{align*}
\hfill $\Box$

Denote by $B_n$ the barbell graph of order $n$ which is obtained from the graph $K_{\lfloor n/2 \rfloor}+K_{\lceil n/2 \rceil}$ by adding one edge.
Note that $B_n=\overline{G_0}.$

{\bf Theorem 10.} {\it Let $G$ be a connected $[4,2]$-graph of order $n\ge 8.$ Then $e(G)\ge \lfloor (n-1)^2/4 \rfloor+1,$ with equality if and only if $G=B_n.$
}

{\bf Proof.}
Since $G$ is a $[4,2]$-graph, $\alpha(G)\le 3.$
If $\alpha(G)=1,$ then $G=K_n$ and $e(G)=n(n-1)/2>\lfloor (n-1)^2/4 \rfloor+1.$
Suppose $\alpha(G)=2.$ Then $\overline{G}$ is triangle-free. If $\overline{G}=K_{\lfloor n/2 \rfloor,\lceil n/2\rceil}$ or $\overline{G}=K_{n/2-1,n/2+1},$ then $G$ is not connected, a contradiction. By Lemma 6, we have $e(\overline{G})\le \lfloor n^2/4 \rfloor-1.$ Combining this inequality with Lemma 8, we have
\[
e(G)=\frac{n(n-1)}{2}-e(\overline{G})\ge \frac{n(n-1)}{2}-\left(\left\lfloor\frac{n^2}{4}\right\rfloor-1\right)=\left\lfloor\frac{(n-1)^2}{4}\right\rfloor +1,
\]
with equality if and only if $\overline{G}=G_0,$ i.e., $G=B_n.$

Suppose $\alpha(G)=3.$
Let $A$ be an independent set of cardinality three. Since $G$ is a $[4,2]$-graph,
for every $v\in V(G)\setminus A,$ we have $|N_A(v)|\ge 2.$
The condition $n\ge 8$ implies that $|G-A|\ge 5.$ By Theorem 9, $e(G-A)\ge \lfloor (n-4)^2/4\rfloor.$ Since $n\ge 8,$ we have
\begin{align*}
e(G)\ge e(G-A)+2(n-3)
&\ge \left\lfloor \frac{(n-4)^2}{4}\right\rfloor+2(n-3)\\
&>\left\lfloor\frac{(n-1)^2}{4}\right\rfloor+1.
\end{align*}
\hfill $\Box$

We denote by $B_n^+$ the $2$-connected graph of order $n$ obtained from the graph $K_{\lfloor n/2 \rfloor}+K_{\lceil n/2 \rceil}$ by adding two edges.

{\bf Theorem 11.} {\it Let $G$ be a $2$-connected $[4,2]$-graph of order $n\ge 10.$ Then $e(G)\ge \lfloor (n-1)^2/4\rfloor+2,$ with equality if and only if $G=B_n^+.$
}

{\bf Proof.} Since $G$ is a $[4,2]$-graph, $\alpha(G)\le 3.$
If $\alpha(G)=1,$ then $G=K_n$ and hence $e(G)=n(n-1)/2>\lfloor (n-1)^2/4\rfloor+2.$
Suppose $\alpha(G)=2.$ Then $\overline{G}$ is triangle-free. If $\overline{G}$ is isomorphic to one of
$$
K_{\lfloor n/2\rfloor,\lceil n/2\rceil},\,\, K_{n/2-1,n/2+1},\,\, K_{n/2-1,n/2+1}^-,\,\, K_{(n-3)/2,(n+3)/2},\,\, G_0,\,\, G_1,\,\, G_2,
$$
then $G$ is not $2$-connected. By Lemmas 6 and 8, we have
\[
e(G)=\frac{n(n-1)}{2}-e(\overline{G})\ge \frac{n(n-1)}{2}-\left(\left\lfloor\frac{n^2}{4}\right\rfloor-2\right)=\left\lfloor\frac{(n-1)^2}{4}\right\rfloor +2,
\]
with equality if and only if $\overline{G}=G_3,$ i.e., $G=B_n^+.$

Suppose $\alpha(G)=3.$ Let $A$ be an independent set of cardinality three. Since $G$ is a $[4,2]$-graph,
for every $v\in V(G)\setminus A,$ we have $|N_A(v)|\ge 2.$ By Theorem 9, we have $e(G-A)\ge \lfloor (n-4)^2/4\rfloor.$ Since $n\ge 10,$ we have
\begin{align*}
e(G)\ge e(G-A)+2(n-3)
&\ge \left\lfloor \frac{(n-4)^2}{4}\right\rfloor+2(n-3)\\
&>\left\lfloor\frac{(n-1)^2}{4}\right\rfloor+2.
\end{align*}
\hfill $\Box$

A graph that contains exactly one Hamilton cycle is called a {\it uniquely hamiltonian graph}.
An interesting unsolved conjecture of Sheehan [10] states that every 4-regular graph is not a uniquely hamiltonian graph.
We have the following result on $[4,2]$-graphs.

{\bf Theorem 12.} {\it
Every $2$-connected $[4,2]$-graph of order at least eight is not uniquely hamiltonian.
}

{\bf Proof.} To the contrary, assume that $G$ is a $2$-connected $[4,2]$-graph of order $n\ge 8$ with a unique Hamilton cycle. By Theorem 5, $G$ contains a cycle with length $n-1,$ say $C.$ Let $v$ be the unique vertex outside $C.$ Specify an orientation of $C.$ We distinguish two cases.

{\bf Case 1.} $\deg_G(v)\ge 3.$

We first consider the case where there exist two neighbors of $v$ that are consecutive on $C.$
Let $v_1,v_2$ be two neighbors of $v$ with $v_1v_2\in E(C).$ Assume that $v_2=v_1^+.$ Let $v_3$ be a neighbor of $v$ distinct from $v_1$ and $v_2.$
If $v_3=v_2^+$ or $v_3=v_1^-,$ then $G$ contains two Hamilton cycles, a contradiction. Thus $d_C(v_2,v_3)\ge 2$ and $d_C(v_1,v_3)\ge 2.$
Since $n\ge 8,$ we have either $d_C(v_2,v_3)\ge 3$ or $d_C(v_1,v_3)\ge 3.$
Without loss of generality, assume that $d_C(v_1,v_3)\ge 3.$ If one of $vv_1^-,vv_2^+,vv_3^+,v_2^+v_3^+$ belongs to $E(G),$ then $G$ contains two Hamilton cycles, a contradiction.
Since $G[\{v,v_1^-,v_2^+,v_3^+\}]$ contains at least two edges, we have $v_1^-v_2^+\in E(G).$
If $d_C(v_2,v_3)=2,$ then $G$ contains two Hamilton cycles, a contradiction. Hence $d_C(v_2,v_3)\ge 3.$
Since $e(\{v,v_2^+,v_3^+,v_2^{+2}\})\ge 2,$ we deduce that one of  $vv_2^+,vv_2^{+2},vv_3^+,v_3^+v_2^+,v_2^{+2}v_3^+$ belongs to $E(G),$ implying that $G$ contains two Hamilton cycles, a contradiction.

Next, we assume that any two neighbors of $v$ are not consecutive on $C.$ Let $v_1,v_2,v_3$ be three neighbors of $v.$ Since $e(\{v_1^+,v_2^+,v_3^+,v\})\ge 2,$ we deduce that two of $$
vv_1^+,\,\, vv_2^+,\,\, vv_3^+,\,\, v_1^+v_2^+,\,\, v_1^+v_3^+,\,\, v_2^+v_3^+
$$
belong to $E(G),$  implying that $G$ contains two Hamilton cycles, a contradiction.

{\bf Case 2.} $\deg_G(v)=2.$

Let $N_G(v)=\{v_1,v_2\}.$ Suppose $d_C(v_1,v_2)=1$ and assume that $v_2=v_1^+.$
Since $G[\{v,v_2^+,v_2^{+3},v_2^{+5}\}]$ contains at least two edges, we have either $v_2^+v_2^{+3}\in E(G)$ or $v_2^{+3}v_2^{+5}\in E(G).$
In the former case, there exists a new $(n-1)$-cycle avoiding $v_2^{+2}.$
If $\deg_G(v_2^{+2})\ge 3,$ then using an argument similar to the one in Case 1, we deduce that $G$ is not uniquely hamiltonian, a contradiction.
If $\deg_G(v_2^{+2})=2,$ then $G[\{v,v_2^{+2},v_2^{+4},v_2^{+5}\}]$ contains exactly one edge $v_2^{+4}v_2^{+5},$ a contradiction.
In the latter case, there exists a new $(n-1)$-cycle avoiding $v_2^{+4}.$
If $\deg_G(v_2^{+4})\ge 3,$ then as in Case 1, we deduce that $G$ is not uniquely hamiltonian, a contradiction.
If $\deg_G(v_2^{+4})=2,$ then $G[\{v,v_2^{+4},v_2^+,v_2^{+2}\}]$ contains exactly one edge $v_2^+v_2^{+2},$ a contradiction.

Suppose $d_C(v_1,v_2)=2$ and assume that $v_2=v_1^{+2}.$
Then there exists a new $(n-1)$-cycle $v_1\overleftarrow{C}v_2vv_1.$
If $\deg_G(v_1^+)\ge 3,$ then as in Case 1, we also deduce a contradiction.
If $\deg_G(v_1^+)=2,$ then $G[\{v,v_1^+,v_2^+,v_2^{+2}\}]$ contains exactly one edge $v_2^+v_2^{+2},$ a contradiction.

Suppose $d_C(v_1,v_2)\ge 3.$
We assert that $v_1^-v_2^{-2}\notin E(G).$
Otherwise, $vv_1\overrightarrow{C}v_2^{-2}v_1^-\overleftarrow{C}v_2v$ is a new $(n-1)$-cycle.
If $\deg_G(v_2^-)\ge 3,$ then as in Case 1, we would deduce a contradiction.
If $\deg_G(v_2^-)=2,$ then $G[\{v,v_2^-,v_2^+,v_2^{+2}\}]$ contains exactly one edge $v_2^+v_2^{+2},$ a contradiction.
Thus $v_1^-v_2^{-2}\notin E(G).$
Since $e(\{v,v_1^-,v_2^-,v_2^{-2}\})\ge 2,$ we have $v_1^-v_2^-\in E(G).$
Similarly, we have $v_1^+v_2^+\in E(G),$ and hence $G$ contains two Hamilton cycles, a contradiction. \hfill$\Box$

The graph in Figure 2 is a $2$-connected $[4,2]$-graph of order seven with a unique Hamilton cycle.
\begin{figure}[h]
\centering
\includegraphics[width=0.4\textwidth]{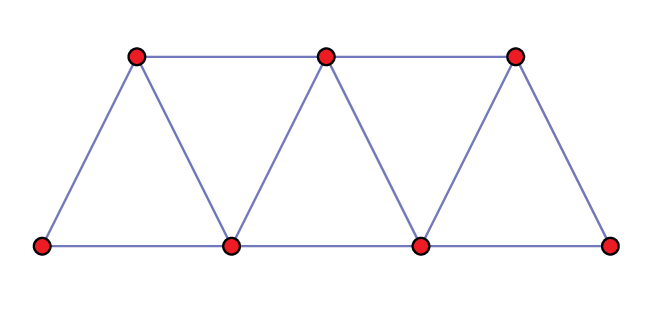}
\caption{A uniquely hamiltonian $[4,2]$-graph of order seven}
\end{figure}
Thus the order requirement in Theorem 12 is necessary.

\section{Concluding remarks}

In this section we mention several unsolved problems posed by the second author who has discussed them with colleagues in private conversations.
The following result was conjectured by Erd\H{o}s and proved by H\"{a}ggkvist, Faudree and Schelp [7]. It generalizes works by Dirac, Ore and Bondy.

{\bf Theorem 13} [7]. {\it If $G$ is a nonbipartite hamiltonian graph of order $n$ with size at least $\lfloor (n-1)^2/4\rfloor +2,$ then $G$ is pancyclic.}

If $n\ge 3$ is an odd integer, we denote by $BT(n)$ the graph obtained by identifying an edge of $K_{(n-1)/2, (n-1)/2}$ with an edge in $C_3.$
The graph $BT(7)$ is depicted in Figure 3.
\begin{figure}[h]
\centering
\includegraphics[width=0.4\textwidth]{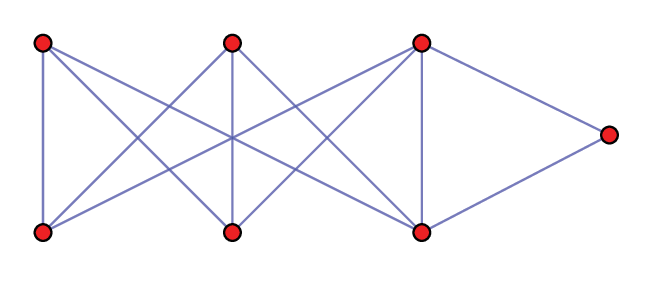}
\caption{The graph $BT(7)$}
\end{figure}

Perhaps Theorem 13  might be strengthened as follows.

{\bf Conjecture 2} (X. Zhan, July 2025). If $G$ is a nonbipartite hamiltonian graph of order $n\ge 7$ with size at least $\lfloor (n-1)^2/4\rfloor +2$ other than
$BT(n)$ when $n$ is odd, then $G$ contains a pancyclic edge.

In view of Theorem 2 and Theorem 11, Conjecture 2, if true, would imply Theorem 5.

The graph in Figure 4 is a vertex-pancyclic graph of order $12$ that contains no pancyclic edge.
\begin{figure}[h]
\centering
\includegraphics[width=0.33\textwidth]{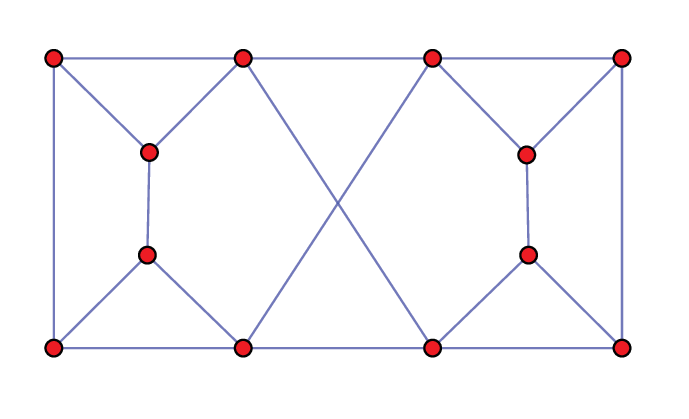}
\caption{A vertex-pancyclic graph that contains no pancyclic edge}
\end{figure}

It is natural to think of the following problem.

{\bf Problem 3} (X. Zhan, October 2025). Characterize the vertex-pancyclic graphs that contain no pancyclic edge.

Clearly, if a vertex-pancyclic graph contains no pancyclic edge, then its minimum degree is at least three.
The following question is a subproblem of Problem 3.

{\bf Question 4} (X. Zhan, October 2025). Does there exist a positive integer $k$ such that every $k$-connected vertex-pancyclic graph contains a pancyclic edge?

\vskip 5mm
{\bf Acknowledgement.} This research  was supported by the NSFC grant 12271170 and Science and Technology Commission of Shanghai Municipality grant 22DZ2229014.

\end{document}